%% file: golden.tex
\documentclass[11pt]{article}

\usepackage{amssymb,amsthm,amsmath}
\usepackage{pstricks,pstricks-add}

\def\Q{\mathord{\mathbb Q}}
\def\R{\mathord{\mathbb R}}
\def\C{\mathord{\mathbb C}}

\def\chr{\mathop{\text{chr}}\nolimits}
\def\arcchr{\mathop{\text{arcchr}}\nolimits}

\newtheorem{theo}{Theorem}
\newtheorem{prop}[theo]{Proposition}
\newtheorem{lemma}[theo]{Lemma}

\title{The Golden Angle is not Constructible}
\author{Pedro J.\ Freitas}
\date{December 14, 2020}

\begin{document}

\maketitle

The golden mean is usually defined with relation to a line segment. A line segment is said to be divided according to the golden ratio if it is decomposed into two segments, with lengths $a>b$, satisfying
\begin{equation}
\frac{a+b}{a} = \frac{a}{b}
\label{goldensection}
\end{equation}
If this happens, the value of these two ratios is the {\em golden number}, $\varphi=(1+\sqrt{5})/2\approx 1.618034$. So, 
$$a=\frac 1 \varphi (a+b)\qquad b = (a+b)-a= \left(1-\frac1\varphi \right) (a+b)$$

The same can be done with a circle, instead of a line segment. A circle is divided into two arcs $\alpha$ and $\beta$ according to the golden ratio if they satisfy equation \eqref{goldensection}, which leads to
\begin{equation}
\alpha=\frac 1 \varphi 2\pi \qquad \beta = \left(1-\frac1\varphi \right) 2\pi \label{angles}
\end{equation}
The smaller angle $\beta$ is called the {\em golden angle} and has some connections to plant growth and phyllotaxis, see \cite[ch. 14]{Th}. Its measure in degrees, approximated to two decimal points, is $137.51^{\rm o}$.
\begin{center}
\input{angulo_ouro.tex}

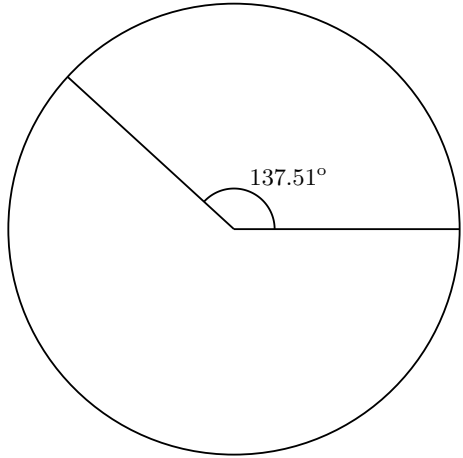
Figure 1. The golden angle
\end{center}
In this note we prove that the golden angle is not constructible with straightedge and compass, by proving that its sine and cosine are transcendental numbers. Since all constructible numbers have to be algebraic, this is enough to prove what we want. 

We recall that the algebraic numbers form a subfield of $\C$, which is closed for taking $n$-th roots. 

\begin{lemma} Given $x\in \R$, we have that $\sin x$ and $\cos x$ are either both algebraic or both transcendental. 

Moreover, the number $e^{ix}$ is transcendental iff either $\cos x$ or $\sin x$ is transcendental (in which case, both are). 
\begin{proof} If both $\sin x$ and $\cos x$ are transcendental, then the first statement is true. If one of them is algebraic, say $\sin x$, then $\cos x = \pm \sqrt{1-\sin^2x}$ is also algebraic. 

For the second statement, if $z=e^{ix}$ is algebraic, then $\cos x = (z+\bar z)/2$ is also algebraic, and similarly for $\sin x$. Conversely, if both $\cos x$ and $\sin x$ are algebraic, then $e^{ix} = \cos x + i \sin x$ is algebraic. 
\end{proof}
\end{lemma} 

We now make use of the Gelfond-Schneider theorem, which is a very powerful tool to generate transcendental numbers (one could also use the Lindemann–Weierstrass theorem). 

\begin{theo}[Gelfond-Schneider] Let $a$ and $b$ be algebraic numbers, such that $a\notin \{0,1\}$ and $b\in \C\setminus \Q$. Then $a^b$ is a transcendental number. 
\end{theo}

See \cite[p.\ 868]{book} as a reference. This theorem solved part of Hilbert's seventh problem, on the irrationality and transcendence of certain numbers

Now consider the angles $\alpha$ and $\beta$ in equation \eqref{angles}. We wish to prove that the golden angle has transcendental sine and cosine.

\begin{prop} 
The golden angle has transcendental sine and cosine, and therefore it is not constructible with straightedge and compass. 
\begin{proof}
Since $\beta=2\pi -\alpha$, it has the same cosine as $\alpha$, and symmetric sine, we can prove that $\alpha=2\pi /\varphi$ has transcendental sine and cosine. For this we prove that $z=e^{i\alpha}=e^{2i\pi /\varphi}$ is transcendental, which is equivalent, according to the lemma. 

If $z$ were algebraic, then, according to the Gelfond-Schneider theorem, $z^{\varphi} = e^{2\pi i} =1$ would be transcendental, which is false. Therefore, $e^{2\pi /\varphi}$ is transcendental.
\end{proof}
\end{prop} 

This proves the non-constructibility of the golden angle. Nevertheless, it is possible to achieve very good approximations, using straightedge and compass. Portuguese artist Almada Negreiros (1893--1970) devoted several years to finding geometric constructions which related to his own analysis of artistic artefacts (see \cite{FC} for more information). One of his discoveries was precisely an approximate construction for the golden angle, presented in  figure 2, based on the regular pentagram, which is constructible. 
\begin{center}
\input{ps-C25-42.tex}

Figure 2. An approximate construction for the golden angle
\end{center}
In this drawing, point $C$ is obtained through an arc of circle centred at $A$. The golden angle is approximated by circle arc $BC$. To compute its exact measure, one only needs to notice three  facts---the two first ones are known properties of the regular pentagon. 
\begin{itemize}
\item The segments marked $a$ and $b$ are proportioned according to the golden number: $a/b=\varphi$. 
\item Length $a$ coincides with the side of the pentagon, that is, it is the chord of $2\pi/5$. 
\item Arc $AC$ has chord $b$.
\end{itemize}

To compute the value of arc $BC$, it is useful to use the {\em chord} as a trigonometric function of the angle. 
\begin{center}
\input{corda1.tex}

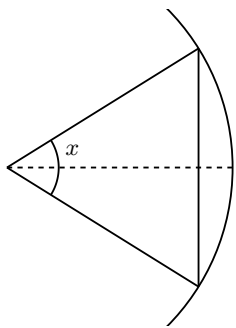
Figure 3. The chord function
\end{center}

Figure 3 helps to deduce its expression, as well as that of its inverse function: 
$$\chr(x) = 2\sin \frac{x}2 \qquad \arcchr(x) = 2\arcsin \frac{x}2$$

Using the three facts mentioned above, we get: 
$$a=2\sin \frac{\pi}5\qquad b = \frac2\varphi \sin \frac{\pi}5$$
$$AC=\arcchr b = 2\arcsin\left( \frac1\varphi \sin \frac{\pi}5\right)$$
From this, we get that the measure of arc $BC$, in degrees, rounded to two decimal values, is $137.40^{\rm o}$. This represents an error of 0.08\% with respect to the golden angle. 

\bibliographystyle{plain}

\end{document}

%% file: angulo_ouro.tex
\psset{xunit=3.0cm,yunit=3.0cm,linewidth=0.7pt,algebraic=true}
\fontsize{9pt}{9pt}
\begin{pspicture*}(-1.3,-1.1)(1.3,1.1)
\pscircle(0.,0.){3.}
\parametricplot{0.0}{2.4000022544174024}{0.180123646705*cos(t)+0.|0.180123646705*sin(t)+0.}
\psline(-0.737395238315,0.675461518156)(0.,0.)
\psline(0.,0.)(1.,0.)
\rput[bl](0.0690567989007,0.193778568393){$137.51^{\rm o}$}
\end{pspicture*}

%% file: ps-C25-42.tex
\psset{xunit=0.15cm,yunit=0.15cm,algebraic=true,linewidth=0.7pt}
\begin{pspicture*}(-25.,-23.)(25.,23.)
\fontsize{9pt}{9pt}
\pscircle(0.,0.){3.}
\psline(-20.,0.)(16.1803398875,-11.7557050458)
\psline(16.1803398875,-11.7557050458)(-6.1803398875,19.0211303259)
\psline(-6.1803398875,19.0211303259)(-6.1803398875,-19.0211303259)
\psline(-6.1803398875,-19.0211303259)(16.1803398875,11.7557050458)
\psline(16.1803398875,11.7557050458)(-20.,0.)
\parametricplot{0.3141592653589793}{1.1990196937026876}{14.5308505601*cos(t)-20 |14.5308505601*sin(t)}
\psline(-20.,0.)(20.,0.)
\rput[tl](-22.5,1){$A$}
\rput[tl](21,1){$B$}
\rput[tl](-17.3,15.3){$C$}
\rput[tl](-13,4.5){$b$}
\rput[tl](3,10){$a$}

\begin{scriptsize}
\psdots[dotstyle=*](-6.1803398875,4.4902797658)
\psdots[dotstyle=*](20.,0.)
\psdots[dotstyle=*](-20.,0.)
\psdots[dotstyle=*](-14.72135955,13.5381524958)
\psdots[dotstyle=*](16.1803398875,11.7557050458)
\end{scriptsize}
\end{pspicture*}

%% file: corda1.tex
\psset{xunit=3.0cm,yunit=3.0cm,linewidth=0.7pt,algebraic=true}
\begin{pspicture*}(2.9,-0.7)(4.1,0.7)
\fontsize{9pt}{9pt}
\pscircle(3.,0.){3.}
\psline(3.84905471491,0.52830492245)(3.84905471491,-0.52830492245)
\psline(3.,0.)(3.84905471491,0.52830492245)
\psline(3.,0.)(3.84905471491,-0.52830492245)
\psline[linestyle=dashed,dash=2pt 2pt](3.,0.)(4.,0.)
\parametricplot{-0.556602893344074}{0.556602893344074}{0.229380048231*cos(t)+3.|0.229380048231*sin(t)+0.}
\rput[tl](3.26,0.1){$x$}
\end{pspicture*}